\documentstyle[11pt,amssymb,amstex]{amsart}
\numberwithin{equation}{section}

\def\bc{{\mathbb C}}

\def\bn{{\mathbb N}}

\def\bq{{\mathbb Q}}

\def\bz{{\mathbb Z}}

\def\g{\gamma}

\def\e{\varepsilon}
\def\l{\lambda}

\def\r{\rho}
\def\s{\sigma}

\newtheorem{thm}{Theorem}[section]
\newtheorem{lem}[thm]{Lemma}
\newtheorem{cor}[thm]{Corollary}
\newtheorem{prop}[thm]{Proposition}

\begin{document}

\title[$p$-adic dynamical system]
{On the chaotic behavior of a generalized logistic $p$-adic
dynamical system}
\author{Farrukh Mukhamedov}
\address{Farrukh Mukhamedov\\
Departamento de Fisica\\
Universidade de Aveiro\\
Campus Universitario de Santiago\\
3810-193 Aveiro, Portugal}\email{{\tt far75m@@yandex.ru},{\tt
farruh@@fis.ua.pt}}

\author{Jos\'{e} F.F. Mendes}
\address{Jos\'{e} F.F. Mendes\\
Departamento de Fisica \\
Universidade de Aveiro\\
Campus Universit\'{a}rio de Santiago\\
3810-193 Aveiro, Portugal} \email{{\tt jfmendes@@fis.ua.pt}}

\begin{abstract}

In the paper we describe basin of attraction $p$-adic dynamical
system $G(x)=(ax)^2(x+1)$. Moreover, we also describe the Siegel
discs of the system, since the structure of the orbits of the
system is related to the geometry of the $p$-adic Siegel discs.

\vskip 0.3cm \noindent {\bf Mathematics Subject Classification}:
37E99, 37B25, 54H20, 12J12.\\
{\bf Key words}:  Basin of attraction, Siegel disc, $p$-adic
dynamics.
\end{abstract}

\maketitle
\small
\section{Introduction}

Applications of $p$-adic numbers in $p$-adic mathematical physics
\cite{ADFV, FW, MP, V1, V2}, quantum mechanics \cite{Kh1} and many
others \cite{Kh2, VVZ} stimulated increasing interest in the study
of $p$-adic dynamical systems.  Note that the $p$-adic numbers
were first introduced by the German mathematician K.Hensel. During
a century after their discovery they were considered mainly
objects of pure mathematics. Starting from 1980's various models
described in the language of $p$-adic analysis have been actively
studied.

On the other hand, the study of $p$-adic dynamical systems arises
in Diophantine geometry in the constructions of canonical heights,
used for counting rational points on algebraic vertices over a
number field, as in \cite{CS}. In \cite{Kh4, TVW} The $p$-adic
field have arisen in physics in the theory of superstrings,
promoting questions about their dynamics. Also some applications
of $p$-adic dynamical systems to some biological, physical systems
were proposed in \cite{ABKO, AKK1, AKK2, DGKS, Kh5, KhN}. In
\cite{BM},\cite{Li} dynamical systems (not only monomial) over
finite field extensions of the $p$-adic numbers were considered.
Other studies of non-Archimedean dynamics in the neighborhood of a
periodic and of the counting of periodic points over global fields
using local fields appeared in \cite{HY,H, L,LP, P}.  Certain
rational $p$-adic dynamical systems were investigated in
\cite{KM},\cite{M},\cite{MR1}, which appear from problems of
$p$-adic Gibbs measures \cite{KMR,MR2,MR3,MRM}. Note that in
\cite{RL,B3,Bz} a general theory of $p$-adic rational dynamical
systems over complex $p$-adic filed $\bc_p$ has been developed. In
\cite{B1, B2} the Fatou set of a rational function defined over
some finite extension of $\bq_p$ has been studied. Besides, an
analogue of Sullivan's no wandering domains' theorem for $p$-adic
rational functions, which have no wild recurrent Julia critical
points, was proved.

The most studied discrete $p$-adic dynamical systems (iterations
of maps) are the so called monomial systems. In
\cite{AKTS},\cite{Kh22} the behavior  of a $p$-adic dynamical
system $f(x)=x^n$ in the fields of $p$-adic numbers $\bq_p$ and
$\bc_p$ was investigated. In \cite{KhN} perturbated monomial
dynamical systems defined by functions $f_q(x) = x^n+q(x)$, where
the perturbation $q(x)$ is a polynomial whose coefficients have
small $p$-adic absolute value, have been studied. It was
investigated the connection between monomial and perturbated
monomial systems.  Formulas for the number of cycles of a specific
length to a given system and the total number of cycles of such
dynamical systems were provided. These investigations show that
the study of perturbated dynamical systems is important. Even for
a quadratic function $f(x)=x^2+c$, $c\in\bq_p$ its chaotic
behavior is complicated (see \cite{TVW, AKK2,Sh}). In \cite{Sh,
DSV} the Fatou and Julia sets of such a $p$-adic dynamical system
were found. Certain ergodic and mixing properties of monomial and
perturbated dynamical systems have been considered in
\cite{A},\cite{GKL}.

The aim of this paper is to investigate the asymptotic behavior of
a nonlinear $p$-adic dynamical system, especially a generalized
$p$-adic logistic map $G(x)=(ax)^2(x+1)$. Note that the logistic
map $f (x)=Cx(1+x)$  and generalized logistic maps are well known
in the literature and it is of great importance in the study of
dynamical systems (see \cite{AK,D,JS}). Much is known about the
behavior of the dynamics of the orbits of a $p$-adic analog of the
logistic map (see \cite{Sh},\cite{TVW}). On the other hand, our
dynamical system is also a perturbated cubic dynamical system,
since it can be reduced to the form $f(x)=x^3+ax^2$. In the paper
we will consider all possible cases of the perturbated term $ax^2$
with respect to the parameter $a$. Note that globally attracting
sets play an important role in dynamics, restricting the
asymptotic behavior to certain regions of the phase space.
However, descriptions of the global attractor can be difficult as
it may contain complicated chaotic dynamics. Therefore, in the
paper we will investigate the basin of attraction of such a
dynamical system. Moreover, we also describe the Siegel discs of
the system, since the structure of the orbits of the system is
related to the geometry of the $p$-adic Siegel discs (see
\cite{AV}).

\section{Preliminaries}

\subsection{$p$-adic numbers}

Let $\bq$ be the field of rational numbers. Throughout the paper
$p$ will be a fixed prime number. Every rational number $x\neq 0$
can be represented in the form $x=p^r\frac{n}{m}$, where
$r,n\in\bz$, $m$ is a positive integer and $p,n,m$ are relatively
prime. The $p$-adic norm of $x$ is given by $|x|_p=p^{-r}$ and
$|0|_p=0$. This norm satisfies so called the strong triangle
inequality
$$
|x+y|_p\leq\max\{|x|_p,|y|_p\}.
$$

From this inequality one can infer that
\begin{eqnarray}\label{inq1}
&&\textrm{if} \ |x|_p\neq |y|_p, \ \ \textrm{then} \ \
|x-y|_p=\max\{|x|_p,|y|_p\}\\[3mm] \label{inq2}&&\textrm{if} \
|x|_p=|y|_p, \ \textrm{then} \ |x-y|_p\leq |2x|_p.
\end{eqnarray}

 This is a ultrametricity of the norm. The completion
of $\bq$ with respect to the $p$-adic norm defines the $p$-adic
field which is denoted by $\bq_p$. Note that any $p$-adic number
$x\neq 0$ can be uniquely represented in the canonical series:
\begin{equation}\label{rep}
x=p^{\g(x)}(x_0+x_1p+x_2p^2+...) ,
\end{equation}
where $\g=\g(x)\in\bz$ and $x_j$ are integers, $0\leq x_j\leq
p-1$, $x_0>0$, $j=0,1,2,...$ (see more detail \cite{Ko,R, S}).
Observe that in this case $|x|_p=p^{-\g(x)}$.

We recall that an integer $a\in \bz$ is called {\it a quadratic
residue modulo $p$} if the equation $x^2\equiv a(\textrm{mod
$p$})$ has a solution $x\in \bz$.

\begin{lem}\label{sq}\cite{S},\cite{VVZ} In order that the equation
\[
x^2=a, \ \ 0\neq a=p^{\g(a)}(a_0+a_1p+...), \ \ 0\leq a_j\leq p-1,
\ a_0>0
\]
has a solution $x\in \bq_p$, it is necessary and sufficient that
the following conditions are satisfied:
\begin{itemize}
\item[(i)] $\g(a)$ is even; \item[(ii)] $a_0$ is a quadratic
residue modulo $p$ if $p\neq 2$, if $p=2$ besides $a_1=a_2=0$.
\end{itemize}
\end{lem}

For any $a\in\bq_p$ and $r>0$ denote
$$
\bar B_r(a)=\{x\in\bq_p : |x-a|_p\leq r\},\ \ B_r(a)=\{x\in\bq_p :
|x-a|_p< r\},
$$
$$
S_r(a)=\{x\in\bq_p : |x-a|_p= r\}.
$$

A function $f:B_r(a)\to\bq_p$ is said to be {\it analytic} if it
can be represented by
$$
f(x)=\sum_{n=0}^{\infty}f_n(x-a)^n, \ \ \ f_n\in \bq_p,
$$ which converges uniformly on the ball $B_r(a)$.

Note the basics of $p$-adic analysis, $p$-adic mathematical
physics are explained in \cite{Ko,R,S, VVZ}

\subsection{Dynamical systems in $\bq_p$}

In this subsection we recall some standard terminology of the
theory of dynamical systems (see for example
\cite{PJS},\cite{KhN}).

Consider a dynamical system $(f,B)$ in $\bq_p$, where $f: x\in
B\to f(x)\in B$ is an analytic function and $B=B_r(a)$ or $\bq_p$.
Denote $x^{(n)}=f^n(x^{(0)})$, where $x^0\in B$ and
$f^n(x)=\underbrace{f\circ\dots\circ f(x)}_n$.
 If $f(x^{(0)})=x^{(0)}$ then $x^{(0)}$
is called a {\it fixed point}. A fixed point $x^{(0)}$ is called
an {\it attractor} if there exists a neighborhood
$U(x^{(0)})(\subset B)$ of $x^{(0)}$ such that for all points
$y\in U(x^{(0)})$ it holds
$\lim\limits_{n\to\infty}y^{(n)}=x^{(0)}$, where $y^{(n)}=f^n(y)$.
If $x^{(0)}$ is an attractor then its {\it basin of attraction} is
$$
A(x^{(0)})=\{y\in \bq_p :\ y^{(n)}\to x^{(0)}, \ n\to\infty\}.
$$
A fixed point $x^{(0)}$ is called {\it repeller} if there  exists
a neighborhood $U(x^{(0)})$ of $x^{(0)}$ such that
$|f(x)-x^{(0)}|_p>|x-x^{(0)}|_p$ for $x\in U(x^{(0)})$, $x\neq
x^{(0)}$. For a fixed point  $x^{(0)}$ of a function $f(x)$ a ball
$B_r(x^{(0)})$ (contained in $B$) is said to be a {\it Siegel
disc} if each sphere $S_{\r}(x^{(0)})$, $\r<r$ is an invariant
sphere of $f(x)$, i.e. if $x\in S_{\r}(x^{(0)})$ then all iterated
points $x^{(n)}\in S_{\r}(x^{(0)})$ for all $n=1,2\dots$. The
union of all Siegel discs with the center at $x^{(0)}$ is said to
{\it a maximum Siegel disc} and is denoted by $SI(x^{(0)})$.

{\bf Remark 2.1.} In non-Archimedean geometry, a center of a disc
is nothing but a point which belongs to the disc, therefore, in
principle, different fixed points may have the same Siegel disc.

Let $x^{(0)}$ be a fixed point of an analytic function $f(x)$. Set
$$
\l=\frac{d}{dx}f(x^{(0)}).
$$

The point $x^{(0)}$ is called {\it attractive} if $0\leq |\l|_p<1$, {\it indifferent} if
$|\l|_p=1$, and {\it repelling} if $|\l|_p>1$.

\section{A generalized logistic map and its fixed points}

Our main interest is a $p$-adic generalized map, which is defined
by
$$
G(x)=(ax)^2(x+1),
$$
where $x,a\in\bq_p$. Using a simple conjugacy $h(x)=ax$ we can
reduce $G$ to the following form $f=h\circ G\circ h^{-1}$, i.e.
\begin{equation}\label{func}
f(x)=x^3+a x^{2}, \ \ \ \ a\in\bq_p.
\end{equation}

Henceforth, we will deal with the function $f$.  Direct checking
shows that the fixed points of the function \eqref{func} are the
following ones
\begin{equation}\label{fix}
x_{1}=0 \ \ \textrm{and} \ \ x_{2,3}=\frac{-a\pm\sqrt{a^2+4}}{2}.
\end{equation}

Here $x_{2,3}$ are the solutions of
\begin{equation}\label{eq1}
x^2+ax-1=0.
\end{equation}

Note that these fixed points are formal, because, basically in
$\bq_q$ the square root does not always exist. A full
investigation of a behavior of the dynamics of the function needs
the existence of the fixed points $x_{2,3}$. Therefore, we have to
verify when $\sqrt{a^2+4}$ does exist.

\subsection{Existence of fixed points}

In this subsection, basically we are going to use Lemma \ref{sq}
to show the existence of $\sqrt{a^2+4}$. Therefore, we consider
several distinct cases with respect to the parameter $a$ and the
prime $p$.

{\bf Case $|a|_p<1$}

In this case we can write $a$ as follows
\begin{equation}\label{form1}
a=p^{k}\e, \ \ k\geq 1, \ \ |\e|_p=1.
\end{equation}
Let us represent $a^2+4$ in the canonical form (see \eqref{rep})
\begin{equation}\label{can1}
a^2+4=p^{\g}(a_0+a_1p+a_2p^2+...) ,
\end{equation}

Now first assume that $p=3$. Then from \eqref{form1} we get
$$
a^2+4=1+3+3^{2k}\e^{2}.
$$
Hence from \eqref{can1} we find that $\g=0$ and $a_0=1$. According
to Lemma \ref{sq} we have to solve the equation $x^2\equiv
1(\textrm{mod $3$})$. One can see that it has a solution $x=3N+1$,
$N\in \bz$. Therefore, in this setting $\sqrt{a^2+4}$ exists.

Let $p\geq 5$. Then we have
$$
a^2+4=4+p^{2k}\e^{2}.
$$
Whence from \eqref{can1} one sees that $\g=0$ and $a_0=4$. The
equation $x^2\equiv 4(\textrm{mod $p$})$ has solution $x=pN+2$,
$N\in \bz$, hence in this setting $\sqrt{a^2+4}$ also exists.

Note that the case $p=2$ is usually pathological, therefore it
should be considered in more detail.

Now let $p=2$.  Then from \eqref{form1} we find
\begin{equation}\label{form11}
a^2+4=2^2(1+2^{2(k-1)}\e^{2}),
\end{equation}
here $|\e|_2=1$. Using Lemma \ref{sq} and $|\e|_2=1$ one gets that
\begin{equation}\label{form12}
\e^2=1+2^m\e_1,
\end{equation}
for some $m\geq 3$, $|\e_1|_2=1$. Now substituting \eqref{form12}
to \eqref{form11} we obtain
\begin{equation}\label{form13}
a^2+4=2^2(1+2^{2(k-1)}+2^{2(k-1)+m}\e_1).
\end{equation}

If $k\geq 3$ then $\g=2$ and $a_0=1, a_1=a_2=0$ in terms of
\eqref{can1}. Therefore, according to Lemma \ref{sq} we infer that
$\sqrt{a^2+4}$ exists.

If $k=2$, then from \eqref{form13} one yields
\begin{equation}\label{form14}
a^2+4=2^2(1+2^{2}+2^{2+m}\e_1).
\end{equation}
Hence, from \eqref{form14} we conclude that $a_1=1$ and therefore
Lemma \ref{sq} implies that $\sqrt{a^2+4}$ does not exist.

Finally, if $k=1$, then we find that
\begin{equation}\label{form15}
a^2+4=2^3(1+2^{m-1}\e_1),
\end{equation}
which with $m\geq 3$ implies that $\g=3$, and again using Lemma
\ref{sq} we conclude that $\sqrt{a^2+4}$ does not exist.

Thus we can formulate

\begin{prop}\label{11} Let $|a|_p<1$, then the expression $\sqrt{a^2+4}$ exists in $\bq_p$ if and
only if either $p\geq 3$ or $p=2$ and $|a|_p\leq 1/p^3$.
\end{prop}

{\bf Case $|a|_p>1$}

In this case we can write that
\begin{equation}\label{form2}
a=p^{-k}\e, \ \ k\geq 1, \ \ |\e|_p=1.
\end{equation}
Then from \eqref{form2} we get
$$
a^2+4=p^{-2k}(\e^{2}+4p^{2k}),
$$
which with Lemma \ref{sq} implies that $\sqrt{a^2+4}$ exists.

\begin{prop}\label{12} Let $|a|_p>1$, then the expression $\sqrt{a^2+4}$ exists in $\bq_p$.
\end{prop}

{\bf Case $|a|_p=1$}

This case is more complicated than others. So we will consider
several subcases with respect to $p$. Before going to details from
$|a|_p=1$ and \eqref{rep} we infer that $a$ can be represented as
follows
\begin{equation}\label{rep1}
a=a_0+a_1p+a_2p^2+\cdots...
\end{equation}
here $a_0\neq 0$.

First assume that $p=2$. Then taking into account Lemma \ref{sq}
we can write
\begin{equation}\label{form3}
a^2=1+b_32^3+b_42^4+\cdots...
\end{equation}
whence one gets
$$
a^2+4=1+2^2+b_32^3+b_42^4+\cdots
$$
But again according to Lemma \ref{sq} we conclude that
$\sqrt{a^2+4}$ does not exist.

Let $p=3$. Then in this case we have
\begin{equation}\label{form31}
a^2=1+b_13+b_23^2+\cdots
\end{equation}
hence
$$
a^2+4=2+(b_1+1)3+b_23^2+\cdots
$$
Consequently, by means of Lemma \ref{sq} one gets that
$\sqrt{a^2+4}$ does not exist, since the equation $x^2\equiv
2(\textrm{mod $3$})$ has not solution in $\bz$.

Let $p\geq 5$. Then from \eqref{rep1} we find
\begin{equation}\label{form311}
a^2+4=a_0^2+4+2a_0a_1p+(a_1^2+2a_0a_2)p^2+\cdots
\end{equation}
Let $a_0^2+4\equiv \!\!\!\!\!\!/ \ 0(\textrm{mod $p$})$ then Lemma
\ref{sq} implies that $\sqrt{a^2+4}$ exists if and only if the
following relation holds
\begin{equation}\label{form32}
\frac{x^2-a^2_0-4}{p}\in\bz
\end{equation}
for some $x\in\bz$.

Let $a_0^2+4\equiv 0(\textrm{mod $p$})$, then we have $a_0^2+4=Np$
for some $N\in\{1,\dots,p-1\}$. From \eqref{form311} one finds
that
\begin{equation}\label{form33}
a^2+4=p(N+2a_0a_1+(a_1^2+2a_0a_2)p+\cdots).
\end{equation}
Taking into account Lemma \ref{sq} and \eqref{form33} we can
formulate the following condition: if
\begin{equation}\label{form322}
\left\{
\begin{array}{lll}
a_0^2+4\equiv 0(\textrm{mod $p$}),\\[2mm]
\frac{a_0^2+4}{p}+2a_0a_1\equiv 0(\textrm{mod $p$}),\\[3mm]
a_1^2+2a_0a_2\equiv \!\!\!\!\!\!/ \ 0(\textrm{mod $p$}),
\end{array}
\right.
\end{equation}
then $\sqrt{a^2+4}$ exists.

For example, if $p=5$, then using the last condition we find that
for the element
$$
a=1+2\cdot 5+2\cdot 5^2+\cdots
$$
$\sqrt{a^2+4}$ exists.

Thus we have the following

\begin{prop}\label{13} Let $|a|_p=1$, then the expression $\sqrt{a^2+4}$ does not exist in
$\bq_p$ when $p=2$ or $p=3$. Let $p\geq 5$. If $|a^2+4|_p=1$, then
$\sqrt{a^2+4}$ exists in $\bq_p$ if and only if \eqref{form32} is
valid for some $x\in\bz$. If $|a^2+4|_p<1$ and \eqref{form322}
holds, then $\sqrt{a^2+4}$ exists.
\end{prop}

\section{Behavior of the fixed points}

In this section we are going to calculate norms of the fixed
points and their behavior.

Let us first note that the derivative of $f$ is
\begin{equation}\label{der1}
f'(x)=3x^2+2ax.
\end{equation}

From this we immediately conclude that the fixed point $x_1$ is
attractive. Therefore, furthermore we will deal with $x_{2,3}$.

Using \eqref{eq1} we find
\begin{equation}\label{viet}
x_2+x_3=-a, \ \ \ \ x_2x_3=-1.
\end{equation}
and
\begin{equation}\label{der2}
f'(x_\s)=3-ax_\s, \qquad \s=2,3.
\end{equation}

{\bf Case $|a|_p<1$}

Let $p\geq 3$, then from \eqref{fix} and \eqref{viet} one finds
that $|x_\s|_p=1$.

Let $p=2$, then Proposition \ref{11} implies that $a=p^k\e$ for
some $k\geq 3$ with $|\e|_p=1$. From this and taking into account
\eqref{fix} we have
\begin{eqnarray*}
|x_\s|_p&=&p|a\pm\sqrt{a^2+4}|_p\\
&=&p|p^k\e\pm\sqrt{p^{2k}\e^2+p^2}|_p\\
&=&|p^{k-1}\e\pm\sqrt{p^{2(k-1)}\e^2+1}|_p=1,
\end{eqnarray*}
since $|p^{2(k-1)}\e^2+1|_p=1$ and $k\geq 3$.

Now let us compute $|f'(x_\s)|_p$. Let $p\neq 3$, then using
\eqref{der2} one gets
\begin{equation*}
|f'(x_\s)|_p=|3-ax_\s|_p=1, \ \ \s=2,3.
\end{equation*}
This means that the fixed points  $x_{\s}$,($\s=2,3$) are
indifferent.

Now let $p=3$, then we easily obtain that $|f'(x_\s)|_p<1$, $
\s=2,3$, which implies that the fixed points are attractive.

We summarize the considered case by the following

\begin{lem}\label{bih1}
Let $|a|_p<1$. The fixed point $x_1$ is attractive. If $p\neq 3$
then the fixed points $x_{2,3}$ are indifferent. If $p=3$ then
$x_{2,3}$ are attractive.
\end{lem}

{\bf Case $|a|_p=1$}

In this case according to Proposition \ref{13} we have to consider
$p\geq 5$. If $|x_2|_p<1$ then \eqref{viet} implies that
$|x_3|_p>1$. This with \eqref{inq1} yields that $|x_2+x_3|_p>1$,
which contradicts to the first equality of \eqref{viet}. Hence
$|x_\s|_p=1$, $\s=2,3$.

Now by means of \eqref{der2} and \eqref{viet} one finds that
\begin{eqnarray}\label{der3}
|f'(x_2)f'(x_3)|_p&=&|(3-ax_2)(3-ax_3)|_p\nonumber \\
&=&|9-3(x_2+x_3)a+a^2x_2x_3|\nonumber \\
&=&|9+2a^2|_p \end{eqnarray} Analogously,
\begin{eqnarray}\label{der4}
|f'(x_2)+f'(x_3)|_p=|6+a^2|_p. \end{eqnarray}

These two equalities \eqref{der3},\eqref{der4} imply that
\begin{eqnarray*}
|f'(x_2)f'(x_3)|_p\leq 1, \ \ \ |f'(x_2)+f'(x_3)|_p\leq 1.
\end{eqnarray*}

Let us prove the following
\begin{lem}\label{atr} Let $|a|_p=1$, then the inequalities
\begin{eqnarray}\label{der5}
|f'(x_2)f'(x_3)|_p< 1, \ \ \ |f'(x_2)+f'(x_3)|_p<1
\end{eqnarray}
are not valid in the same time.
\end{lem}
\begin{pf} Let us assume that \eqref{der5} is valid. According to \eqref{der3},\eqref{der4}
from \eqref{der5} we obtain $|9+2a^2|_p<1$ and $|6+a^2|_p<1$. The
last ones equivalent to the following equations
\begin{eqnarray*}
9+2x^2\equiv 0(\textrm{mod $p$}), \ \ \ 6+x^2\equiv 0(\textrm{mod
$p$}).
\end{eqnarray*}
From these relations we infer that $p=3$, which is impossible
thanks to Proposition \ref{13}.
\end{pf}

This Lemma implies that the both fixed points can not be
simultaneously attractive.  Now let us provide some examples for
the occurrence of the other cases.\\

{\sc Example 4.1.}  Let $p=5$. Then according to Proposition
\ref{13} we infer that $a_0=1$ or $a_0=4$ in representation
\eqref{rep1}. So we have $a^2=1+p\e_1$ for some $|\e_1|_p=1$.
Consequently, \eqref{der3} and \eqref{der4} imply that
\begin{equation*}
|f'(x_2)f'(x_3)|_p=1, \ \ \ |f'(x_2)+f'(x_3)|_p=1,
\end{equation*}
which means that $|f'(x_2)|_p=1$, $|f'(x_3)|_p=1$, i.e. both fixed
points are indifferent.

If $p=7$ we also similarly can get analogous  result as the
previous one, i.e. $x_\s$ is indifferent.\\

{\sc Example 4.2.} Let $p=11$ and $a=4$. Then Proposition \ref{13}
implies that $\sqrt{a^2+4}$ exists. On the other hand, from the
equalities \eqref{der3}, \eqref{der4} we find that
\begin{equation*}
|f'(x_2)f'(x_3)|_p=1, \ \ \ |f'(x_2)+f'(x_3)|_p<1,
\end{equation*}
which yields that $|f'(x_2)|_p=1$, $|f'(x_3)|_p=1$, i.e. both
fixed points are indifferent. Note also that in this case from
\eqref{der4} we get $|a^2+4|_p=1$.\\

{\sc Example 4.3.} Now let $p=11$ and $a=1$. Then Proposition
\ref{13} implies that $\sqrt{a^2+4}$ exists. Consequently, from
\eqref{der3} and \eqref{der4} we infer that
\begin{equation*}
|f'(x_2)f'(x_3)|_p<1, \ \ \ |f'(x_2)+f'(x_3)|_p=1,
\end{equation*}
which implies that either $|f'(x_2)|_p<1$, $|f'(x_3)|_p=1$ or
$|f'(x_2)|_p=1$, $|f'(x_3)|_p<1$. Without loss of generality we
can assume that $|f'(x_2)|_p<1$, $|f'(x_3)|_p=1$, which means that
$x_2$ is attractive and $x_3$ is indifferent. Note that in this
case from \eqref{der3} we have
\begin{equation}\label{aa}
|2(a^2+4)+1|_p=|2a^2+9|<1 \ \ \ \Rightarrow \ \ |a^2+4|_p=1.
\end{equation}\\

Thus we can formulate the following

\begin{lem}\label{bih2}
Let $|a|_p=1$. The fixed point $x_1$ is attractive. For the other
fixed points there are the following possibilities:
\begin{itemize}
    \item[1.] both fixed points $x_2$ and $x_3$ are indifferent;
    \item[2.] the fixed point $x_2$ is attractive and $x_3$ is indifferent,
respectively.
    \item[3.] the fixed point $x_3$ is attractive and $x_2$ is
indifferent, respectively.\\
\end{itemize}
\end{lem}

{\bf Case $|a|_p>1$}

In this case from \eqref{viet} one gets that $|x_2+x_3|_p=|a|_p$,
$|x_2|_p|x_3|_p=1$. These imply that either $|x_2|_p>1$ or
$|x_3|_p>1$. Without lose of generality we may assume that
$|x_2|_p>1$, which means that $|x_3|_p<1$.  From the condition
$|a|_p>1$ we find that $|x_2|_p=|a|_p$ and $|x_3|_p=1/|a|_p$.

From \eqref{der2} we obtain
\begin{equation}\label{der6}
|f'(x_2)|_p=|a|_p|x_2|_p=|a|_p^2>1
\end{equation}
which means that the point $x_2$ is repelling.

From \eqref{der1} with $|x_3|_p=1/|a|_p$ one gets that
$|f'(x_3)|_p=1$, hence $x_3$ is an indifferent point.

\begin{lem}\label{bih3}
Let $|a|_p>1$. The fixed point $x_1$ is attractive. For the other
fixed points there are the following possibilities:
\begin{itemize}
       \item[1.] the fixed point $x_2$ is repelling and $x_3$ is indifferent,
respectively.
    \item[2.] the fixed point $x_3$ is repelling and $x_2$ is
indifferent, respectively.\\
\end{itemize}
\end{lem}

\section{Attractors and Siegel discs}

In the previous section we have established behavior of the fixed
points of the dynamical system. Using those results, in this
section we are going to describe the size of attractors and Siegel
discs of the system.

Before going to details let us formulate certain useful auxiliary
facts.

Let us assume that $x^{(0)}$ is  a fixed point of $f$. Then $f$
can be represented as follows
\begin{equation}\label{tay}
f(x)=f(x^{(0)})+f'(x^{(0)})(x-x^{(0)})+\frac{f''(x^{(0)})}{2}(x-x^{(0)})^2+\frac{f'''(x^{(0)})}{6}(x-x^{(0)})^3.
\end{equation}
From the above equality putting $\g=x-x_0$ we obtain
\begin{equation}\label{dif1}
|f(x)-f(x^{(0)})|_p=|\g|_p\bigg|f'(x^{(0)})+\frac{f''(x^{(0)})}{2}\g+\frac{f'''(x^{(0)})}{6}\g^2\bigg|_p.
\end{equation}

\begin{lem}\label{att1}
Let $x^{(0)}$ be a fixed point of the function $f$ given by
\eqref{func}. If for $\g=x-x^{(0)}$ the following inequality holds
\begin{equation*}
\max\bigg\{|3x^{(0)}+a|_p|\g|_p,|\g|^2_p\bigg\}<|f'(x^{(0)})|_p
\end{equation*}
then
\begin{equation}\label{dif2}
|f(x)-f(x^{(0)})|_p=|\g|_p|f'(x^{(0)})|_p.
\end{equation}
\end{lem}

The proof immediately comes from \eqref{dif1} and the following
ones
\begin{equation*}
f''(x)=6x+2a, \ \ \ \ f'''(x)=6.
\end{equation*}
From Lemma \ref{att1} we get

\begin{cor}\label{att2}\cite{AKTS} Let $x^{(0)}$ be a fixed point of the function $f$ given by
\eqref{func}. The following assertions hold:
\begin{itemize}
    \item[(i)] if $x^{(0)}$ is an attractive point of $f$, then it is an attractor
of the dynamical system. If $r>0$ satisfies the inequality
\begin{equation}\label{att3}
\max\bigg\{|3x^{(0)}+a|_pr,r^2\bigg\}<1
\end{equation}
then $B_r(x^{(0)})\subset A(x^{(0)})$;
    \item[(ii)] if $x^{(0)}$ is an indifferent point of $f$ then it is the
center of a Siegel disc. If $r$ satisfies the inequality
\eqref{att3} then $B_r(x^{(0)})\subset SI(x^{(0)})$;
    \item[(iii)] if $x^{(0)}$ is a repelling point of $f$ then $x^{(0)}$ is a
repeller of the dynamical system.
\end{itemize}
\end{cor}

Now as in the previous section we consider several distinct cases
with respect to the parameter $a$.

{\bf Case $|a|_p>1$}

In the previous section point out that the fixed point $x_1=0$ is
attractive (see Lemma \ref{bih3}), therefore let us first
investigate $A(x_1)$. To do it, denote
$$
r_k=\frac{1}{|a|^k_p}, \ \ k\geq 0.
$$

Now consider several steps along the description of $A(x_1)$.

{\tt (I)}. From Corollary \ref{att2} and \eqref{att3} we find that
$B_{r_1}(0)\subset A(x_1)$. Now take $x\in S_{r_1}(0)$, i.e.
$|x|=r_1$. Then one gets
$$
|f(x)|_p=|x|^2_p|x+a|_p=|x|^2_p|a|_p=r_1,
$$
whence we infer that $|f^{(n)}(x)|_p=r_1$ for all $n\in\bn$. This
means that $x\not\in A(x_1)$, hence $A(x_1)\cap
S_{r_1}(0)=\emptyset$. Moreover, we have $f(S_{r_1}(0))\subset
S_{r_1}(0)$.

In the sequel we will assume that
\begin{equation}\label{ddd}
\sqrt{|a|_p}\notin\{p^k, \ k\in\bn\}. \end{equation}
Denote
\begin{equation*}
A(\infty)=\{x\in\bq_p: |f^{(n)}(x)|_p\to\infty \ \ \textrm{as} \ \
n\to\infty\}.
\end{equation*} It is evident that $A(x_1)\cap
A(\infty)=\emptyset$.

{\tt (II)}. Let us take $x\in S_{r}(0)$ with $r>|a|_p$. Then we
have
$$
|f(x)|_p=|x|^2_p|x+a|_p=|x|^2_p|x|_p=|x|^3_p,
$$
which means that $x\in A(\infty)$, i.e. $S_{r}(0)\subset
A(\infty)$ for all $r>|a|_p$.

{\tt (III)}. Now assume that $x\in S_r(0)$ with $r\in
(r_1,r_0)\cup (r_0,|a|_p)$. Then we have $f(S_r(0))\subset
S_{r^2|a|_p}(0)$. If $r\in (r_0,|a|_p)$ then $r^2|a|_p>|a|_p$,
hence we have $S_r(0)\subset A(\infty)$. If $r\in (r_1,r_0)$ then
according to our assumption \eqref{ddd} we have
$r^{2^n}|a|_p^{1+2+\cdots+2^{n-1}}\neq 1$ for every $n\in\bn$,
hence there is $n_0\in\bn$ such that $f^{(n_0)}(S_r(0))\subset
A(\infty)$, from this one concludes that $S_r(0)\subset
A(\infty)$. Consequently, we have $S_{r}(0)\subset A(\infty)$ for
all $r\in (r_1,r_0)\cup (r_0,|a|_p)$.

{\tt (IV)}. If $x\in S_{r_0}(0)$, then one gets
$f(S_{r_0}(0))\subset S_{|a|_p}(0)$.

{\tt (V)}. Therefore, we have to consider $x\in S_{|a|_p}(0)$.
From \eqref{func} we can write
\begin{equation}\label{ff}
|f(x)|_p=|a|^2_p|x+a|_p.
\end{equation}
From this we conclude that we have to investigate behavior of
$|x+a|_p (\leq |a|_p)$. It is clear the following decomposition
\begin{equation}\label{dec}
S_{|a|_p}(0)=\bigcup_{r=0}^{|a|_p}S_{r}(-a).
\end{equation}

{\tt (VI)}. Now if $|x+a|_p<r_3$ then from the last equality we
get $|f(x)|_p<r_1$, this yields that $x\in A(x_1)$. Hence
$B_{r_3}(-a)\subset A(x_1)$. Moreover, taking into account (I) we
have $f(S_{r_3}(-a))\subset S_{r_1}(0)$.

{\tt (VII)}. If $x\in S_{r_1}(-a)$ then from \eqref{ff} we find
that $f(x)\in S_{|a|_p}(0)$.

{\tt (VIII)} If $x\in S_{r_2}(-a)$ then again using \eqref{ff} one
gets that $f(x)\in S_{r_0}(0)$. This with (IV) implies that
$f^{(2)}(S_{r_1}(-a))\subset S_{|a|_p}(0)$.

{\tt (IX)} If $x\in S_{r}(-a)$ with $r\in (r_3,r_2)\cup
(r_2,r_1)\cup (r_1,|a|_p]$. Then from \eqref{ff} we obtain that
$f(x)\in S_{\rho}(0)$, $\rho\in (r_1,r_0)\cup (r_0,|a|_p)\cup
(|a|_p,|a|^3_p]$. Hence thanks to (II) and (III) we infer that
$f(x)\in A(\infty)$.

Let us introduce some more notations. Given sets $A,B\subset\bq_p$
put
\begin{eqnarray}\label{TD1}
T_{f,A,B}(x)&=&\min\{k\in\bn: \ f^{(k)}(x)\in B\}, \ \ x\in A, \\
\label{TD2} D[A,B]&=&\{x\in A: \ T_{f,A,B}(x)<\infty\}.
\end{eqnarray}

Taking into account (II)-(IX) and \eqref{TD1}-\eqref{TD2} we can
define $D[S_{r_0}(0)\cup S_{|a|_p}(0),B_{r_3}(-a)]$, which is non
empty, since from (VI) one sees that $B_{r_3}(-a)\subset
D[S_{r_0}(0)\cup S_{|a|_p}(0),B_{r_3}(-a)]$. Thus taking into
account (I)-(IX) we have the following equality
$A(x_1)=B_{r_1}(0)\cup D[S_{r_0}(0)\cup
S_{|a|_p}(0),B_{r_3}(-a)]$.\\

Now turn to the other fixed points. According to Lemma \ref{bih3}
without loss of generality we may assume that $x_2$ is repelling
and $x_3$ is indifferent.  In this case we know  that
$|x_2|_p=|a|_p$ and $|x_3|_p=r_1$. Now Corollary \ref{att2} with
\eqref{att3} yields that $B_{r_1}(x_3)\subset SI(x_3)$. It is
clear that $0\notin SI(x_3)$, therefore $SI(x_3)=B_{r_1}(x_3)$.
From (I) we infer that $SI(x_3)\subset S_{r_1}(0)$.

Thus we have proved the following

\begin{thm} Let $|a|_p>1$ and \eqref{ddd} be satisfied. Then
$A(x_1)=B_{r_1}(0)\cup D[S_{r_0}(0)\cup
S_{|a|_p}(0),B_{r_3}(-a)]$. For the indifferent fixed point
$x_{3}$ we have $SI(x_3)=B_{r_1}(x_3)$.
\end{thm}

{\bf Case $|a|_p<1$}

From Lemma \ref{bih1} we know that $x_1$ is attractive. So
according to Corollary \ref{att2} we immediately find that
$B_1(0)\subset A(x_1)$. Take $x\in S_1(0)$ then
$|f(x)|_p=|x|^2_p|x+a|_p=|x|^3=1$, hence $|f^{(n)}(x)|_p=1$ for
all $n\in\bn$. This means that $x\not\in A(x_1)$, hence
$A(x_1)=B_{1}(0)$.

Now turn to the fixed points $x_2$ and $x_3$. According to Lemma
\ref{bih1} we consider  two possible situations $p\neq 3$ and
$p=3$.

First assume $p\neq 3$. In this case $x_2$ and $x_3$ are
indifferent, so Corollary \ref{att2} again implies that
$B_{1}(x_\s)\subset SI(x_\s)$, $\s=2,3$.

Let us take $x\in S_r(x_\s)$, $r\geq 1$, then put $\g=x-x_\s$. It
is clear that $|\g|_p=r$. By means of \eqref{dif1} and \eqref{eq1}
we find
\begin{eqnarray}\label{dif3}
|f(x)-f(x_\s)|_p&=&|\g|_p|3x_\s^2+2ax_\s+(3x_\s+a)\g+\g^2|_p
\nonumber
\\
&=&r|\g^2+3x_\s\g+3+a(\g-x_\s)|_p.
\end{eqnarray}

If $r>1$ then from \eqref{dif3} we easily obtain that
$$
|f(x)-f(x_\s)|_p=|\g|_p^3
$$
since $|\g^2+3x_\s\g+3|_p=|\g|_p^2$, $|a(\g-x_\s)|_p=|a|_p|\g|_p$.
This implies that $SI(x_\s)\subset \bar B_1(x_\s)$.

\begin{lem}\label{sig}
Let $|a|_p<1$ and $p\neq 3$. The equality $SI(x_\s)=\bar
B_1(x_\s)$ holds if and only if the equality holds
\begin{equation}\label{sig1}
|\g^2+3x_\s\g+3|=1
\end{equation}
for every $\g\in S_1(0)$.
\end{lem}

\begin{pf} If \eqref{sig1} is satisfied for all
$\g\in S_1(0)$ then from \eqref{dif3} we infer that
$f(S_1(x_\s))\subset S_1(x_\s)$, since $|a(\g-x_\s)|_p<1$. This
proves the assertion. Now suppose that $SI(x_\s)=\bar B_1(x_\s)$
holds. Assume that \eqref{sig1} is not valid, i.e. there is
$\g_0\in S_1(0)$ such that
\begin{equation}\label{sig2}
|\g_0^2+3x_\s\g_0+3|<1
\end{equation}
The last one with \eqref{dif3} implies that $|f(x_0)-x_\s|_p<1$
for an element $x_0=x_\s+\g_0$. But this contradicts to
$SI(x_\s)=\bar B_1(x_\s)$. This completes the proof.
\end{pf}

From the proof of Lemma \ref{sig} we immediately obtain that if
there is $\g_0\in S_1(0)$ such that \eqref{sig2} is satisfied then
$SI(x_\s)=B_1(x_\s)$. Moreover, we can formulate the following

\begin{lem}\label{sig22} Let $|a|_p<1$ and $p\neq 3$. The following conditions are
equivalent:
\begin{itemize}
    \item[(i)] $SI(x_\s)=B_1(x_\s)$;
    \item[(ii)] there is $\g_0\in S_1(0)$ such that \eqref{sig2} is
    satisfied;
    \item[(iii)] $\sqrt{-3}$ exists in $\bq_p$.
\end{itemize}
\end{lem}

\begin{pf} The implication (i)$\Leftrightarrow$(ii) immediately
follows from the proof of Lemma \ref{sig1}. Consider the
implication (ii)$\Rightarrow$(iii). The condition \eqref{sig2}
according to the Hensel Lemma (see \cite{Ko} yields that the
existence of a solution $z\in\bq_p$ of the following equation
\begin{equation}\label{sig3}
z^2+3x_\s z+3=0
\end{equation}
such that $|z-\g_0|<1$ which implies that $|z|_p=1$.

Now assume that there is a solution $z_1\in\bq_p$ of \eqref{sig3}.
Then from Vieta's formula we conclude the existence of the another
solution $z_2\in\bq_p$ such that
$$
z_1+z_2=-3x_\s, \ \ \ z_1z_2=3.
$$
From these equalities one gets that $|z_1+z_2|_p=1$,
$|z_+z_-|_p=1$ which imply $z_1,z_2\in S_1(0)$. So putting
$\g_0=z_1$ we find \eqref{sig2}.

Let us now analyze when \eqref{sig3} has a solution belonging to
$\bq_p$. We know that a general solution of \eqref{sig3} is given
by
\begin{equation}\label{sig4}
z_{1,2}=\frac{-3x_\s\pm\sqrt{-3-9ax_\s}}{2},
\end{equation}
here we have used \eqref{eq1}. But it belongs to $\bq_p$ if
$\sqrt{-3-9ax_\s}$ exists in $\bq_p$. Since $|-9ax_\s|_p=|a|_p<1$
implies that $-9ax_\s=p^k\e$ for some $k\geq 1$ and $|\e|_p=1$.
Hence, $-3-9ax_\s=-3+p^k\e$. Therefore according to Lemma \ref{sq}
we conclude that $\sqrt{-3-9ax_\s}$ exists if and only if
$\sqrt{-3}$ exists in $\bq_p$. The implication
(iii)$\Rightarrow$(ii) can be proven along the reverse direction
in the previous implication.
\end{pf}

Let us consider some concrete examples when $\sqrt{-3}$ exists
with respect to $p$.

{\sc Example 5.1.}Let $p=2$, then $-3$ can be rewritten as follows
$$
-3=1+2^2+2^3+\cdots
$$
so according to Lemma \ref{sq} we concludes that $\sqrt{-3}$ does
not exist in $\bq_2$.

Analogously reasoning we may establish that when $p=5,11$ we find
that $\sqrt{-3}$ does not exist. If $p=7,13$ then $\sqrt{-3}$
exists.\\

From Lemmas \ref{sig} and \ref{sig22} we conclude that $SI(x_\s)$
is either $B_1(x_\s)$ or $\bar B_1(x_\s)$. The equality
\eqref{fix} yields that
\begin{equation}\label{sol}
|x_2-x_3|_p=|\sqrt{a^2+4}|_p=|2|_p,
\end{equation}
which implies that $SI(x_2)\cap SI(x_3)=\emptyset$ when $p\geq 5$
and $SI(x_2)=SI(x_3)$ when $p=2$, since any point of a ball is its
center.\\

Now consider the case $p=3$. According to Lemma \ref{bih1} we see
that the both fixed points  $x_2$ and $x_3$ are attractive. Taking
into account $|x_\s|_p=1$ and $|a|_p<1$ from Corollary \ref{att2}
one finds that $B_{1}(x_\s)\subset A(x_\s)$, $\s=2,3$. From the
equality \eqref{sol} we have $|x_2-x_3|_p=1$ which implies that
$S_1(x_\s)\nsubseteq A(x_\s)$.

Let us take $x\in S_r(x_\s)$ with $r\geq 1$, then putting
$\g=x-x_\s$ from \eqref{dif3} with
$|3x_\s\g+3+a(\g-x_\s)|_p<|\g|_p$ we get
\begin{eqnarray*}
|f(x)-x_\s|_p&=&|\g||\g^2+3x_\s\g+3+a(\g-x_\s)|_p=r^3,
\end{eqnarray*}
which implies that $f(S_r(x_\s))\subset S_{r^3}(x_\s)$ for every
$r\geq 1$. Hence, in particular, we obtain $f(S_1(x_\s))\subset
S_{1}(x_\s)$.

Consequently we have the following

\begin{thm}\label{sig5} Let $|a|_p<1$. The following assertions hold:
\begin{itemize}
    \item[(i)] Basin of attraction of $x_1$ is $B(x_1)$, i.e. $A(x_1)=B_1(0)$;
    \item[(ii)] Let $p\neq 3$. Then $SI(x_\s)=B_1(x_\s)$ is valid if and only if $\sqrt{-3}$
    exists in $\bq_p$. Otherwise $SI(x_\s)=\bar B_1(x_\s)$ holds.
    \item[(iii)] If $p\geq 5$ then $SI(x_2)\cap SI(x_3)=\emptyset$, if
    $p=2$ then $SI(x_2)=SI(x_3)$.
    \item[(iv)] Let $p=3$. Then  $A(x_\s)=B_1(x_\s)$ ($\s=2,3$).
\end{itemize}
\end{thm}

Note that if we consider our dynamical system over $p$-adic
complex field $\bc_p$ we will obtain different result from the
formulated Theorem, since $\sqrt{-3}$ always exists in $\bc_p$.

{\bf Case $|a|_p=1$}

In this case according to Proposition \ref{13} we have to consider
$p\geq 5$.

Let us first describe the basin of attraction of the fixed point
$x_1=0$. Analogously, reasoning as in the previous cases  we may
find that $B_1(0)\subset A(x_1)$. Now if $x\in S_r(0)$ with $r>1$
one gets that $|f(x)|_p=|x|_p^3=r^3$, which implies that $x\notin
A(x_1)$ and $A(x_1)\subset \bar B_1(0)$.

Suppose that $x\in S_1(0)$. From \eqref{func} we find
\begin{equation}\label{1f}
|f(x)|_p=|x+a|_p
\end{equation}
From this we conclude that we have to investigate behavior of
$|x+a|_p$. It is clear the following decomposition
\begin{equation}\label{1dec}
S_{1}(0)=\bigcup_{r=0}^{1}S_{r}(-a).
\end{equation}

Let $|x+a|_p\leq r$ with $r<1$, then from \eqref{1f} we get
$|f(x)|_p<1$, this yields that $x\in A(x_1)$. Hence
$B_{1}(-a)\subset A(x_1)$.

So we can define the set $D[S_{1}(0),B_{1}(-a)]$, which is not
empty since $B_{1}(-a)\subset D[S_{1}(0),B_{1}(-a)]$.

We are going to show that $S_1(0)\setminus D[S_{1}(0),B_{1}(-a)]$
is open. Now assume that $y\in S_1(0)\setminus
D[S_{1}(0),B_{1}(-a)]$. Then this means that
\begin{equation}\label{0df}
|f^{(n)}(y)+a|_p=1 \ \ \textrm{for all} \ \ n\in\bn.
\end{equation}
For any $r<1$ one can be established that $B_r(y)\subset S_1(0)$.
We will show that $B_r(y)\subset S_1(0)\setminus
D[S_{1}(0),B_{1}(-a)]$, which would be the assertion. Take $x\in
B_r(y)$. To show $x\in S_1(0)\setminus D[S_{1}(0),B_{1}(-a)]$ it
is enough to prove $|f^{(n)}(x)+a|_p=1$ for all $n\in\bn$. To do
end, consider
\begin{eqnarray*}
|f(y)-f(x)|_p&=&|y-x|_p|y^2+xy+x^2-a(y+x)|_p\nonumber\\
&\leq &|y-x|_p\leq r
\end{eqnarray*}
which implies that $f(x)\in B_r(f(y))$. By means of the induction
we find that $f^{(n)}(x)\in B_r(f^{(n)}(y))$ for all $n\in\bn$.
This with \eqref{0df} implies the assertion.

Thus we have the following

\begin{thm} Let $|a|_p=1$. Then the fixed point $x_1$ is attractor and
$A(x_1)=B_{1}(0)\cup D[S_{1}(0),B_{1}(-a)]$.
\end{thm}

Let us turn to the fixed points $x_2$ and $x_3$. Note that
furthermore, we always assume that these fixed points exist, which
in accordance with Section 3 is equivalent to the existence of
$\sqrt{a^2+4}$. Now according to Lemma \ref{bih2} we will consider
three different cases.

{\sc Case 1.} In this case $x_\s$ ($\s=2,3$) is indifferent. By
means of Corollary \ref{att2} and using the same procedure as in
the previous cases we immediately derive that $B_1(x_\s)\subset
SI(x_\s)$.

Let $x\in S_r(x_\s)$ with $r>1$. It then follows from \eqref{dif3}
that
\begin{eqnarray}\label{3dif}
|f(x)-f(x_\s)|_p&=&|\g|_p|\g^2+3x_\s\g+3+a(\g-x_\s)|_p=r^3,
\end{eqnarray}
since $r^2=|\g^2+3x_\s\g+3|_p>|a(\g-x_\s)|_p=r$, where as before
$\g=x-x_\s$. This shows that $f(S_r(x_\s))\subset S_{r^3}(x_\s)$,
from which one concludes that $SI(x_\s)\subset \bar B_1(x_\s)$.

Now let $x\in S_1(x_\s)$. Then from \eqref{3dif} we obtain
\begin{eqnarray}\label{4dif}
|f(x)-f(x_\s)|_p&=&|\g^2+\g(3x_\s+a)+3-ax_\s|_p.
\end{eqnarray}

From this we derive the following

\begin{lem}\label{1sig} Let $|a|_p=1$ and $x_\s$ is an indifferent
fixed point of $f$. The the following assertions equivalent:
\begin{itemize}
    \item[(i)] $SI(x_\s)=\bar B_1(x_\s)$;
    \item[(ii)] the equality
\begin{eqnarray}\label{4eq}
|z^2+z(3x_\s+a)+3-ax_\s|_p=1
\end{eqnarray}
holds for all $z\in S_1(0)$;
\end{itemize}
If $|a^2+4|<1$ then the last condition equivalent to
\begin{itemize}
    \item[(iii)] the equality
\begin{eqnarray}\label{5eq}
|2z^2-az+6+a^2|_p=1
\end{eqnarray}
holds for all $z\in S_1(0)$.
\end{itemize}
\end{lem}

\begin{pf} The implication (i)$\Leftrightarrow$(ii) can be proved
along the same line of the proof of Lemma \ref{sig}. Now assume
that $|a^2+4|<1$. So we have to prove the implication
(ii)$\Leftrightarrow$(iii). From \eqref{fix} one gets
\begin{eqnarray}\label{5dif}
3x_\s+a&=&\frac{-a\pm 3\sqrt{a^2+4}}{2},\\
\label{6dif} 3-ax_\s&=&\frac{6+a^2\mp\sqrt{a^2+4}}{2}.
\end{eqnarray}
Using \eqref{5dif}-\eqref{6dif} from \eqref{4eq} we obtain
\begin{eqnarray}\label{7dif}
|z^2+z(3x_\s+a)+3-ax_\s|_p&=&|2z^2+z(-a\pm
3\sqrt{a^2+4})+6+a^2\mp\sqrt{a^2+4}|_p\nonumber\\
&=&|(2z^2-az+6+a^2)\pm(3z-1)\sqrt{a^2+4}|_p.
\end{eqnarray}

Taking in account our assumption from \eqref{7dif} we conclude
that \eqref{4eq} holds if and only if \eqref{5eq} holds.
\end{pf}

\begin{lem}\label{4sig} Let $|a|_p=1$ and $x_\s$ ($\s=2,3$) is an indifferent
fixed point of $f$. The the following assertions equivalent:
\begin{itemize}
    \item[(i)] $SI(x_\s)=B_1(x_\s)$;
    \item[(ii)] $\sqrt{\frac{a^2-6+a\sqrt{a^2+4}}{2}}$ exists in $\bq_p$.
\end{itemize}
\end{lem}

\begin{pf} From Lemma \ref{1sig} we find that (i) is valid if and
only if there is $z_0\in S_1(0)$ such that
\begin{equation}\label{dd2}
|z_0^2+z_0(3x_\s+a)+3-ax_\s|_p<1.
\end{equation}
This thanks to the Hensel Lemma implies that the existence a
solution $b\in\bq_p$ with $|b-z_0|<1$ of the following equation
\begin{eqnarray}\label{8eq}
z_0^2+z_0(3x_\s+a)+3-ax_\s=0.
\end{eqnarray}

On the other hand, the condition $|3-ax_\s|_p=1$ with Vieta's
formulas implies that if a solution $b_1\in\bq_p$ of the equation
exists, then $|b_1|_p=1$ holds. Putting $z_0=b_1$ we get
\eqref{dd2}. So \eqref{dd2} and \eqref{8eq} are equivalent. Hence,
\eqref{8eq} has solution if and only if
$$
\sqrt{\frac{a^2-6+a\sqrt{a^2+4}}{2}}
$$
exists in $\bq_p$, since every solution of \eqref{8eq} has a form
$$
b_{1,2}=\frac{-(3x_\s+a)\pm
\sqrt{\frac{a^2-6+a\sqrt{a^2+4}}{2}}}{2}
$$
\end{pf}

Let us consider more special cases.

\begin{cor}\label{2sig} Let $|a|_p=1$, $|a^2+4|<1$ and $x_\s$ is an indifferent fixed point of $f$. The
the following assertions equivalent:
\begin{itemize}
    \item[(i)] $SI(x_\s)=B_1(x_\s)$;
    \item[(ii)] there is $z_0\in S_1(0)$ such that
\begin{eqnarray}\label{6eq}
|2z_0^2-az_0+6+a^2|_p<1
\end{eqnarray}
holds;
     \item[(iii)] $\sqrt{-5}$ exists in $\bq_p$ at $p\geq 7$.
\end{itemize}

\end{cor}

\begin{pf} The implication (i)$\Leftrightarrow$(ii) is a direct
consequence of Lemma \ref{1sig}. Consider the implication
(ii)$\Leftrightarrow$(iii). In this case \eqref{6eq} thanks to the
Hensel Lemma implies that the existence a solution $b\in\bq_p$
with $|b-z_0|<1$ of the following equation
\begin{eqnarray}\label{7eq}
2z^2-az+6+a^2=0.
\end{eqnarray}

According to Lemma \ref{atr} we get $|a^2+6|=1$, which with
$|a|_p=1$ and Vieta's formulas implies that equivalence of
\eqref{7eq} and \eqref{6eq}. Hence, the existence of \eqref{7eq}
equivalent to the existence of the squire root of the discriminant
of \eqref{7eq}, namely $\sqrt{-48-7a^2}\in\bq_p$.  Rewrite
$-48-7a^2$ as follows
\begin{equation}\label{dd0}
-48-7a^2=-20-7(a^2+4).
\end{equation}
From $|a^2+4|_p<1$ and the existence of $\sqrt{a^2+4}$ we infer
that $a^2+4=p^{2n}\e$ can be written for some $n\in\bn$ and
$|\e|_p=1$. From \eqref{dd0} we have
\begin{equation}\label{dd}
-48-7a^2=-20-p^{2n}\e.
\end{equation}
Hence, if $p\geq 7$, then  keeping in mind \eqref{dd} and Lemma
\ref{sq} we deduce that $\sqrt{-48-7a^2}\in\bq_p$ if and only if
$\sqrt{-5}$ exists in $\bq_p$. If $p=5$ then from \eqref{dd} we
find that $-48-7a^2=p(-4-p^{2n-1}\e)$ which according to Lemma
\ref{sq} yields that $\sqrt{-48-7a^2}$ does not exist in $\bq_p$.
\end{pf}

{\bf Remark 5.1.} If $p=5$ and $|a^2+4|_p<1$ then from the proof
of the last Lemma and Lemma \ref{1sig} we immediately obtain that
$SI(x_\s)=\bar B_1(x_\s)$. \\

Let us turn to the case $|a^2+4|_p=1$. Note that this case  is a
rather tricky. Therefore, we will provide that some more
sufficient conditions for  fulfilling the equality
$SI(x_\s)=B_1(x_\s)$.

Using \eqref{7dif} we assume that $|3z_0-1|_p<1$ for some $z_0\in
S_1(0)$. Then from the equality
\begin{eqnarray}\label{dd1}
|2z_0^2-az_0+6+a^2|_p=\bigg|\frac{2}{3}(3z_0-1)^2+\frac{4-a}{3}(3z_0-1)+\frac{1}{3}(3a^2-a+20)\bigg|_p\
\end{eqnarray}
we conclude that $|2z_0^2-az_0+6+a^2|_p<1$ if and only if
$|3a^2-a+20|_p<1$. By means of \eqref{7dif} and the last
condition, we can formulate the following

\begin{cor}\label{3sig} Let $|a|_p=1$, $|a^2+4|=1$ and $x_\s$ ($\s=2,3$) is an indifferent fixed point of $f$. If
$|3a^2-a+20|_p<1$ holds, then $SI(x_\s)=B_1(x_\s)$.
\end{cor}

Summarizing the obtained  we formulate

\begin{thm}\label{bih4} Let $|a|_1=1$ and $x_2$,$x_3$  are indifferent fixed points of
$f$. The following assertions hold:
\begin{itemize}
    \item[(i)] If $\sqrt{\frac{a^2-6+a\sqrt{a^2+4}}{2}}$ exists in
    $\bq_p$ then $SI(x_\s)=B_1(x_\s)$. Otherwise $SI(x_\s)=\bar
    B_1(x_\s)$.
    \item[(ii)] Let $|a^2+4|_p<1$. If $\sqrt{-5}$ exists in
    $\bq_p$ at $p>5$, then $SI(x_\s)=B_1(x_\s)$. Otherwise $SI(x_\s)=\bar
    B_1(x_\s)$. Moreover, we have $SI(x_2)=SI(x_3)$.
\end{itemize}
\end{thm}

Note that the last assertion immediately follows from
\eqref{sol}.\\

Note that the case 3 similar to the case 2, therefore we will
consider only case 2.

{\sc Case 2.}  In this setting we will suppose that the fixed
point $x_2$ is attractive and $x_3$ is indifferent, respectively.
Recall that later according to \eqref{aa} means that
$|a^2+4|_p=1$. For the point $x_3$ the Siegel discs would the same
as in the previous case. So we have to investigate only $x_2$.

We can easily  show that $B_1(x_2)\subset A(x_2)$. By means of
\eqref{dif3} we can also establish that $A(x_2)\cap
S_r(x_2)=\emptyset$ for all $r>1$.

Let $x\in S_1(x_2)$. Then from \eqref{4dif} one holds
\begin{eqnarray}\label{1a}
|f(x)-f(x_2)|_p&=&|\g^2+3x_2\g+a\g+3-ax_2|_p.
\end{eqnarray}
Attractivity of the point $x_2$ means that $|3-ax_2|_p<1$
therefore if
\begin{eqnarray}\label{2a}
|\g^2+3x_2\g+a\g|_p&=&|\g+3x_2+a|\nonumber\\
&=&|x+2x_2+a|_p<1
\end{eqnarray}
then from \eqref{1a} we get that $x\in A(x_2)$, i.e.
$B_1(-2x_2-a)\subset A(x_2)$. Here we have used the notation
$\g=x-x_2$.

If $|\g^2+3x_2\g+a\g|_p=1$ then $f(x)\in S_1(x_2)$. So we may
again repeat the above procedure. Hence this leads that we can
define the set $D[S_1(x_2),B_1(-2x_2-a)]$, which is nonempty.
Using the same argument as above cases (i.e. $|a|_p<1$) one can
show that $A(x_2)=B_1(x_2)\cup D[S_1(x_2),B_1(-2x_2-a)]$.

\begin{thm} Let $|a|_1=1$. Assume that $x_2$ is attractive and $x_3$ is indifferent fixed points of
$f$. The following assertions hold:
\begin{itemize}
    \item[(i)] $A(x_2)=B_1(x_2)\cup D[S_1(x_2),B_1(-2x_2-a)]$;
    \item[(ii)] The Siegel disc of $x_3$ would be the same as in Theorem \ref{bih4}.
    \end{itemize}
\end{thm}

%\end{document}

\section*{ Acknowledgements} The first author (F.M.) thanks the FCT (Portugal) grant SFRH/BPD/17419/2004.
Another author (J.F.F.M.) acknowledges projects
POCTI/FAT/46241/2002, POCTI/ MAT/46176/2002 and European research
NEST project DYSONET/ 012911.

\end{document}